\documentclass[reqno,a4paper,final]{amsart}
\usepackage{amssymb,amstext,amsthm,eucal}
\usepackage{bbold}
\usepackage{array}
\usepackage[latin1]{inputenc}
\newcommand{\ecke}{\;_-\!\rule{0.2mm}{0.2cm}\;\;}

\newcommand{\U}{\mathrm{U}}

\newcommand{\RR}{\mathbb{R}}

\newtheorem{THEO}{\bf Theorem}
\newtheorem{DEF}{\bf Definition}

\newtheorem{LM}{\bf Lemma}

\newcommand{\MUNCH}[1]{\relax}
\allowdisplaybreaks[1]
\begin{document}
\begin{sloppypar}
\title{On pseudo-Hermitian Einstein spaces}
\author{Felipe Leitner}
\address{Institut f{\"u}r Mathematik, Universit{\"a}t Stuttgart, Germany}
\email{leitner@mathematik.uni-stuttgart.de}
%\thanks{ }
%\date{\today}

\begin{abstract} 
We describe and construct here pseudo-Hermitian  
structures $\theta$ without torsion (i.e. with transversal symmetry)
whose Webster-Ricci curvature tensor is a constant multiple
of the exterior differential $d\theta$. 
We call these structures pseudo-Hermitian Einstein and our result 
states that 
they all can be derived locally from K{\"a}hler-Einstein metrics.
Moreover, we discuss the corresponding Fefferman metrics of the
pseudo-Hermitian Einstein structures. These Fefferman metrics 
are never Einstein, but they are locally always conformally Einstein.
\end{abstract}

\maketitle
\tableofcontents

%%%%%%%%%%%%%%%%%%%%%%%%%%%%%%%%
\section{Introduction}%%%
%%%%%%%%%%%%%%%%%%%%%%%%%%%%%%%%
\label{ab1}
CR-geometry is a $|2|$-graded parabolic geometry on a smooth manifold
$M^n$. Underlying Weyl-structures are the pseudo-Hermitian forms
$\theta$. CR-geometry is closely related to conformal geometry via
the Fefferman construction. For conformal structures, there is
the notion of being conformally Einstein, that means there
is a Riemannian metric in the conformal class which is Einstein. 
In terms of tractor calculus the conformal Einstein condition can be 
expressed
through the existence of a parallel standard tractor (cf. e.g. 
\cite{Gov04}, 
\cite{Lei05}). The concept
of parallel standard tractors works for CR-geometry as well.
One can define in this case that a pseudo-Hermitian structure 
with parallel standard CR-tractor, whose first 'slot' is a 
constant 
real function, is Einstein.

However, we do not use here tractor calculus to define the Einstein 
condition
for a pseudo-Hermitian structure. Instead, we say a pseudo-Hermitian 
structure $\theta$ is Einstein if and only if its torsion vanishes and the
Webster-Ricci curvature is a constant multiple of the exterior 
differential $d\theta$. The two definitions for pseudo-Hermitian
Einstein spaces coincide. 

As our main result, we will show here a construction principle
for pseudo-Hermitian Einstein spaces. In fact, they are closely related with
K{\"a}hler-Einstein spaces (cf. Theorem \ref{TH1}). And we will explicitly show
that the Fefferman metrics which belong to pseudo-Hermitian Einstein spaces 
admit a local Einstein scale (cf. Theorem \ref{TH2}). 

We will proceed as follows.
In section \ref{ab2} we introduce the notions that we use here for 
pseudo-Hermitian geometry, in particular, Webster curvature. 
In section \ref{ab3} we consider pseudo-Hermitian structures with
transversal symmetry and define the Einstein condition.
In section \ref{ab4} we compare the pseudo-Hermitian geometry of $\theta$
with the Riemannian 
geometry of the induced metric $g_\theta$. In section \ref{ab5}
we derive the natural Riemannian submersion of a transversally symmetric
pseudo-Hermitian space. We will see that the Ricci tensor of the base 
space
of the Riemannian submersion determines the Webster-Ricci tensor
of the transversally symmetric pseudo-Hermitian space. 
In section \ref{ab6} we find the construction principles for pseudo-Hermitian 
Einstein spaces taking off with a K{\"a}hler-Einstein space 
(cf. Theorem \ref{TH1}). 
Finally, in section \ref{ab7} and \ref{ab8} we recall the Fefferman 
construction
and prove explicitly the conformal Einstein condition for those Fefferman metrics
which come from pseudo-Hermitian Einstein structures (cf. Theorem 
\ref{TH2}).

%%%%%%%%%%%%%%%%%%%%%%%%%%%%%%%%
\section{Pseudo-Hermitian structures}%%%
%%%%%%%%%%%%%%%%%%%%%%%%%%%%%%%%
\label{ab2}
We fix here in brief some notations for pseudo-Hermitian structures.
Threreby, we follow mainly the notations of \cite{Bau99}. More material
on pseudo-Hermitian geometry can be found e.g. in \cite{Lee86}, 
\cite{Lee88}, \cite{CS00}
\cite{Cap01} or \cite{CG02}.

With a CR-structure on a smooth manifold $M^n$ of odd
dimension $n=2m+1$ we mean here a pair $(H,J)$, which consists of 
\begin{enumerate}
\item a contact distribution $H$ in $TM$ of codimension $1$ and
\item a complex structure $J$ on $H$, i.e. $J^2=-id|_H$, subject to the
integrability conditions
$[JX,Y]+[X,JY]\in\Gamma(H)$ and 
\[J([JX,Y]+[X,JY])-[JX,JY]+[X,Y]=0\]
for all $X,Y\in\Gamma(H)$.
\end{enumerate} 
The conditions that the distribution $H$ is contact and 
the complex structure $J$ is integrable ensures that $(H,J)$ 
determines a $|2|$-graded parabolic geometry on $M$ (cf. e.g. 
\cite{CS00}).
In particular, the (infinitesimal) automorphism group of $(M,H,J)$
is finite dimensional.   

A nowhere vanishing real $1$-form $\theta\in\Omega(M)$ is called a 
pseudo-Hermitian structure on the CR-manifold $(M,H,J)$ if 
\[\theta|_H\equiv 0\ .\] Then we call the data $(M,H,J,\theta)$ a pseudo-Hermitian 
space. Since the distribution $H$ is contact, the $1$-form $\theta$ is 
necessarily
a contact form. Such a contact form $\theta$ exist on $(M,H,J)$ 
if and only if $M$ is orientable. Furthermore, two pseudo-Hermitian 
structures $\theta$ and $\tilde{\theta}$ on 
$(M,H,J)$ differ only by multiplication with a real nowhere vanishing
function $f\in C^\infty(M)$:
\[\tilde{\theta}=f\cdot \theta\ .\]

We consider now the exterior differential $d\theta$ of a pseudo-Hermitian
structure. 
This $2$-form 
is non-degenerate on $H$, i.e.
\[(d\theta)^m|_H\neq 0\ \] and the $2$-tensor
\[L_\theta(\cdot,\cdot):= d\theta(\cdot,J\cdot)\]
is symmetric and non-degenerate on $H$. If $L_\theta$ is positive definite
the pseudo-Hermitian structure $\theta$ is called strictly pseudoconvex.
In general, the $2$-tensor $L_\theta$ has complex signature $(p,q)$ on $H$
(resp. real signature $(2p,2q)$).
The conditions
\[T\ecke\theta\equiv 1\qquad\mbox{and}\qquad T\ecke d\theta\equiv 0\]
uniquely determine a vector field $T$ on $M$. This $T$ is called 
Reeb vector field. For convenience, we set $J(T)=0$. 

To a pseudo-Hermitian structure $\theta$ on $(M,H,J)$ (with arbitrary
signature for $L_\theta$)
belongs a canonical covariant  
derivative
\[\nabla^W:\Gamma(TM)\longrightarrow\Gamma(T^*M\otimes TM)\ ,\]
which is called the Tanaka-Webster connection. It is uniquely determined 
by the following conditions:
\begin{enumerate}
\item
The connection $\nabla^W$ is metric with respect to the non-degenerate 
symmetric 
$2$-tensor \[g_\theta:=L_\theta+\theta\circ\theta\] on $M$,
i.e. \[\nabla^W g_\theta=0\ ,\] and
\item
its torsion 
$Tor^W(X,Y):=\nabla^W_XY-\nabla^W_YX-[X,Y]$
satisfies 
\[\begin{array}{l}Tor^W(X,Y)=L_\theta(JX,Y)\cdot T\qquad\mbox{for\ all}\ 
X,Y\in\Gamma(H)\quad\mbox{and}\\[3mm]
Tor^W(T,X)=-\frac{1}{2}([T,X]+J[T,JX])\qquad\mbox{for \ all}\ 
X\in\Gamma(H)\ .\end{array}\]
\end{enumerate}
In addition, for this connection it holds  \[
\nabla^W\theta=0\qquad\mbox{and}\qquad\nabla^W\circ J=J\circ\nabla^W\ .\]

The curvature operator of the connection $\nabla^W$ is defined in the 
usual manner:
\[R^{\nabla^W}(X,Y)=[\nabla^W_X,\nabla_Y^W]-\nabla_{[X,Y]}^W\ .\]
The $(4,0)$-curvature tensor $R^W$ is given for $X,Y,Z,V\in TM$ by
\[R^W(X,Y,Z,V):=g_\theta(R^{\nabla^W}(X,Y)Z,V)\ .\]
This curvature tensor has the symmetry properties
\[\begin{array}{l}
R^W(X,Y,Z,V)=-R^W(Y,X,Z,V)=-R^W(X,Y,V,Z),\\[2mm]
R^W(X,Y,JZ,V)=-R^W(X,Y,Z,JV) \ .
\end{array}\]
We have not listed here the Bianchi type identities.
We just note that the Bianchi identities for $R^{\nabla^W}$ do 
not 
(formally) look like 
those
for the Riemannian curvature tensor. We will come back to this point 
later.
 
There is also a notion of Ricci curvature for pseudo-Hermitian
structures. It is called the Webster-Ricci curvature tensor and can be 
defined as follows. Let \[(e_\alpha,Je_\alpha)_{\alpha=1,\ldots ,m}\]
be a local orthonormal frame of $L_\theta$ on $H$ and $\varepsilon_\alpha
:=g_\theta(e_\alpha,e_\alpha)$. 
Then it is defined
\[Ric^W(X,Y):=i\sum_{\alpha=1}^m \varepsilon_\alpha 
R^W(X,Y,e_\alpha,Je_\alpha)\ .\]
The Webster-Ricci curvature is skew-symmetric with values in the 
purely imaginary numbers $i\RR$. And the Webster scalar curvature is
\[scal^W:=i\sum_{\alpha=1}^m \varepsilon_\alpha 
Ric^W(e_\alpha,Je_\alpha)\ .\]
The function $scal^W$ on $(M,H,J,\theta)$ is real.

%%%%%%%%%%%%%%%%%%%%%%%%%%%%%%%%
\section{Transversal symmetry}%%%
%%%%%%%%%%%%%%%%%%%%%%%%%%%%%%%%
\label{ab3}

Let $(M,H,J)$ be a CR-manifold. A vector field $T\neq 0$ is called a 
transversal symmetry of $(H,J)$ if it is not tangential to the subbundle 
$H$ (i.e. it is
transversal to $H$) and if the flow of $T$ consists (at least locally for 
small parameters) of CR-automorphisms, i.e.  the 
distribution $H$ is preserved and $\mathcal{L}_TJ=0$, or 
equivalently 
\[ [T,X]+J[T,JX]=0\qquad \mbox{for\ all}\ X\in\Gamma(H)\ .\]

Now let $\theta$ be a non-degenerate pseudo-Hermitian 
structure on $(M,H,J)$ and let $T$ be the corresponding Reeb
vector field determined by
\[\theta(T)\equiv 1\qquad\mbox{and}\qquad T\ecke d\theta\equiv 0\ .\]      
Obviously, the Reeb vector field to $\theta$ is a transversal 
symmetry of $(H,J)$ if and only if the torsion part $Tor^W(T,X)$ of the 
Tanaka-Webster
connection $\nabla^W$ vanishes for all vector fields $X\in\Gamma(H)$.
Equivalently, it is right to say that $T$ is a transversal symmetry if and
only if $T$ is a Killing vector field for the metric $g_\theta$, i.e.
\[\mathcal{L}_Tg_\theta=0\ .\] This uses the fact that 
\[\mathcal{L}_TJ=0\qquad\mbox{and}\qquad 
\mathcal{L}_T\theta=0\]
for the case when $T$ is a transversal symmetry.

The above observations suggest the following notation. 
We say that a non-degenerate pseudo-Hermitian structure $\theta$ on 
a CR-manifold $(M,H,J)$  is transversally symmetric if its Reeb vector 
field $T$ is a transversal symmetry of $(H,J)$. In short, we say
$\theta$ is a (TSPH)-structure on $(M,H,J)$. 

We extend our notations here further and say that
$\theta$ is a pseudo-Hermitian Einstein structure on $(M,H,J)$ 
if and only if $\theta$ is transversally symmetric and the Webster-Ricci
curvature $Ric^W$ is a constant multiple of $d\theta$, i.e.
\[Ric^W=-i\frac{\ scal^W}{m}\cdot d\theta\qquad\mbox{and}\qquad 
Tor^W(T,X)=0\] for all $ X\in\Gamma(H)$.
In this case 
$(M,H,J,\theta)$ is called a pseudo-Hermitian Einstein space (cf. 
\cite{Lee88}).

%%%%%%%%%%%%%%%%%%%%%%%%%%%%%%%%
\section{Comparision between $\nabla^W$ and 
$\nabla^{g_\theta}$ and 
their curvature tensors}%%%
%%%%%%%%%%%%%%%%%%%%%%%%%%%%%%%%
\label{ab4}

We determine in this section the endomorphism
\[D^\theta:=\nabla^W-\nabla^{g_\theta},\]
where $\nabla^{g_\theta}$ denotes the Levi-Civita connection of 
$g_\theta$,
and derive comparision formulas for the Riemannian and the 
Webster curvature tensors. We will restrict this discussion to 
the transversally symmetric case.  

So let $\theta$ be a (TSPH)-structure on $(M,H,J)$. A straightforward 
calculation shows
that the covariant derivative
\[\nabla^{W}-\frac{1}{2}d\theta\cdot T+
\frac{1}{2}(\theta\otimes J+J\otimes\theta)\]
is metric and has no torsion with respect to $g_\theta$. We conclude
that it is the Levi-Civita connection of $g_\theta$ and we obtain as 
comparision tensor
\[D^\theta:=\nabla^W-\nabla^{g_\theta}=\frac{1}{2}\left(
d\theta\cdot T-(\theta\otimes J+ J\otimes \theta)\right)\ .\]
Another straightforward calculation shows that for any 
$X,Y,Z\in\Gamma(TM)$ it holds
\begin{eqnarray*}R^{\nabla^W}(X,Y)Z&=&R^{g_\theta}(X,Y)Z-\frac{1}{2}
\left(\nabla^{g_\theta}_Zd\theta(X,Y)\right)\cdot T-\frac{1}{2}d\theta(X,Y)\cdot 
J(Z)\\[2mm]
&&+\frac{1}{4}d\theta(Y,Z)\cdot J(X)-\frac{1}{4}d\theta(X,Z)\cdot J(Y)
\\[2mm] 
&&+\frac{1}{4}\theta(Z)\cdot\theta(X)\cdot Y
-\frac{1}{4}\theta(Z)\cdot\theta(Y)\cdot X\ .
\end{eqnarray*}
This is the comparision of the curvature tensors.
The formula immediately proves that (in the
transversally symmetric case!) the Webster curvature $R^{\nabla^W}$
resp. $R^W$ satisfies the first Bianchi identity of the style of a
Riemannian 
curvature tensor, i.e. it holds
\[R^{\nabla^W}(X,Y)Z+R^{\nabla^W}(Y,Z)X+R^{\nabla^W}(Z,X)Y=0\ .\]
This is our main observation 
here. 
%(broad laughter everywhere!).
\begin{LM} Let $\theta$ be a (TSPH)-structure on $(M,H,J)$. Then
the Webster curvature tensor $R^W$ satisfies
\[R^{W}(X,Y,Z,V)+R^{W}(Y,Z,X,V)+R^{W}(Z,X,Y,V)=0\] 
for all $X,Y,Z,V\in TM$. In particular, it holds
\[\begin{array}{l}
R^W(X,Y,Z,V)=R^W(Z,V,X,Y)\qquad\mbox{and}\\[1.5mm]
R^W(X,JY,JZ,V)=R^W(JX,Y,Z,JV)\ .\end{array}\]
\end{LM} 
%(the laughter ebbs slowly away, some people are coughing)

Using the derived symmetry properties of the Webster curvature for
the particular case of transversal symmetry, we obtain the 
following comparision 
between the Riemannian Ricci tensor and the Webster-Ricci tensor. 
Let \[(e_\alpha,Je_\alpha)_{\alpha=1,\ldots, m}=(e_i)_{i=1,\ldots, 2m}\]
denote a local orthonormal frame of $H$ in $TM$. It is
\[Ric^{g_\theta}(X,Y)=R^{g_\theta}(X,T,T,Y)+
\sum_{i=1}^{2m}\varepsilon_iR^{g_\theta}(X,e_i,e_i,Y)
%\\[2mm]
%&=&\sum_{\alpha=1}^m\varepsilon_{\alpha}R^{g_\theta}(X,e_\alpha,e_\alpha,Y)
%+\sum_{\alpha=1}^m\varepsilon_\alpha 
%R^{g_\theta}(X,Je_\alpha,Je_\alpha,Y)+R^{g_\theta}(X,T,T,Y)\\[3mm]
%\end{eqnarray*}
\]
and
\begin{eqnarray*}Ric^W(X,Y)&=&i\sum_\alpha\varepsilon_\alpha
R^W(X,Y,e_\alpha,Je_\alpha)\\[2mm]
&=&i\sum_\alpha\varepsilon_\alpha
R^W(Y,e_\alpha,Je_\alpha,X)+i\sum_\alpha\varepsilon_\alpha
R^W(e_\alpha,X,Je_\alpha,Y)\\[2mm]
&=&i\sum_\alpha\varepsilon_\alpha
R^W(JY,Je_\alpha,Je_\alpha,X)+i\sum_\alpha\varepsilon_\alpha 
R^W(X,e_\alpha,e_\alpha,JY)\\[2mm]
&=&i\sum_i\varepsilon_i R^W(X,e_i,e_i,JY)
\end{eqnarray*}
for all $X,Y\in TM$. 
Moreover, by the comparision formula for the curvature tensors
$R^{g_\theta}$ and $R^W$, we have 
\[\sum_i\varepsilon_i 
R^{\nabla^W}(X,e_i)e_i=Ric^{g_\theta}(X)-R^{g_\theta}(X,T)T+\frac{3}{4}X
%-\frac{1}{2}\sum_i\varepsilon_i
%\left(\nabla_{e_i}^{g_\theta}d\theta(X,e_i)\right)\cdot T
\]
for all $X\in\Gamma(H)$ and 
\[\sum_i\varepsilon_i
R^{\nabla^W}(T,e_i)e_i=\sum_i\varepsilon_i\left(R^{g_\theta}(T,e_i)e_i-\frac{1}{2}
\left(\nabla_{e_i}^{g_\theta}d\theta(T,e_i)\right)\cdot T\right)\ .\]
These formulas combined with 
the fact that 
$R^{g_\theta}(X,T)T=\frac{1}{4}X$ for $X\in H$ result to
\[\begin{array}{l}
Ric^{g_\theta}(X,Y)=iRic^W(X,JY)-\frac{1}{2}g_\theta(X,Y),\\[2mm]
Ric^W(T,X)=0, \qquad Ric^W(T,T)=0\end{array}\]
and 
\[Ric^{g_\theta}(T,X)=0\ ,\qquad
Ric^{g_\theta}(T,T)=\frac{m}{2}g_\theta(T,T)\ ,\]
whereby $X,Y\in H$.

%%%%%%%%%%%%%%%%%%%%%%%%%%%%%%%%
\section{The natural Riemannian submersion of a (TSPH)-structure}%%%
%%%%%%%%%%%%%%%%%%%%%%%%%%%%%%%%
\label{ab5}

We assume here that $\theta$ is a (TSPH)-structure on the CR-manifold 
$(M,H,J)$ of dimension $n=2m+1$. This implies that the Reeb vector field
$T$ to $\theta$ is Killing for the induced metric $g_\theta$.
At least locally, we can factorise through the integral curves of $T$ on 
$M$
and obtain a semi-Riemannian metric $h$ on a factor space,
which has dimension $2m$. We describe this process here in more
detail. In particular, we calculate the relation for the Ricci 
curvatures
of the induced metric $g_\theta$ and the factorised metric $h$.

Let $\theta$ be a (TSPH)-structure on $(M,H,J)$ of signature $(p,q)$. To 
every point 
in $p\in M$ exists a neigborhood (e.g. some small 
ball) $U\subset M$ and a map $\phi_U$ such that 
$\phi_U$ is a diffeomorphism between $U$ and the $\RR^n$, and moreover,
it holds $d\phi_U(T)=\frac{\partial}{\partial x_1}$, that is the first 
standard coordinate vector in $\RR^n$. This implies that there exists 
a smooth submersion
\[\pi_U:U\subset M\to N\subset \RR^{2m}\] such that for all $v\in N$
the inverse image $\pi^{-1}_U(v)$ consists of an integral curve of $T$
through some point in $U$ 
parametrised by an interval in $\RR$. 
Since $T$ is a Killing vector field,
the expression \[h(X,Y):=L_\theta(\pi_{U*}^{-1}X,\pi^{-1}_{U*}Y)\]
is uniquely defined for any $X,Y\in TN$ and gives rise to a smooth metric 
tensor on $N$ of dimension $2m$ of signature $(2p,2q)$. 
Alternatively, we can define \[
h(X,Y)=g_\theta(X^*,Y^*)\ ,\]
where $X^*$ denotes the horizontal lift of the vector $X$ to $M$
with respect to $g_\theta$ and the vertical direction $\RR T$.
In particular, the 
map
\[ \pi_U: (U,g_\theta) \to (N,h)\]
is a smooth Riemannian submersion. The construction is naturally 
derived from $\theta$ only (and some chosen neighborhood $U$).
The distribution $H$ in $TU$ is 
horizontal for this submersion (i.e. orthogonal to the vertical).

For simplicity, we assume now that \[ \pi: (M,g_\theta)\to (N,h)\]
is globally a smooth Riemannian submersion, whereby the inverse images are 
the 
integral curves of the Reeb vector field $T$ to a (TSPH)-structure 
$\theta$ on $M$ with CR-structure $(H,J)$. 
Since the complex structure $J$ acts on $H$ and $T$
is an infinitesimal automorphism of $J$, the complex structure can be 
uniquley projected to a smooth endomorphism on $N$, which we
also denote by $J$ and which satisfies $J^2=-id|_N$. Since $J$ is 
integrable on $H$, the endomorphism $J$ is integrable on $N$ as well, i.e.
$J$ is a complex structure on $N$. In fact, $J$ is a K{\"a}hler structure
on $(N,h)$, i.e. \[
\nabla^{h}J=0\ .\]
The latter fact can be seen with the comparision tensor $D^\theta$. It is 
\begin{eqnarray*}
(\nabla^{g_\theta}_{X^*}J)(Y^*) &=& 
\nabla^{g_\theta}_{X^*}(JY^*)-J\nabla_{X^*}^{g_\theta}Y^*\\[3mm]
&=& 
\nabla^W_{X^*}(JY^*)-(J\nabla^W_{X^*}Y^*)-\frac{1}{2}d\theta(X^*,J(Y^*))\cdot 
T\\
&=& -\frac{1}{2}g_\theta(X^*,Y^*)\cdot T
\end{eqnarray*}
and 
\[Vert_{\pi}\nabla_{X^*}^{g_\theta}(J(Y^*))=-\frac{1}{2}g_\theta(Y^*,X^*)\cdot 
T\ .\]
Together with $\nabla^{h}\circ\pi=\pi\circ\nabla^{g_\theta}$ this implies
$\nabla^{h}J=0$
on $N$. 

Altogether, we know yet that a (TSPH)-space $(M,H,J,\theta)$ gives rise 
(locally) in a natural manner
to a $(2m)$-dimensional K{\"a}hler space $(N,h,J)$. 
We use now the 
well known
formulas for the Ricci tensor of a Riemannian submersion to calculate
$Ric^{h}$ (cf. \cite{ONe66}). 
The application of the standard formulas shows that
\[\begin{array}{l}Ric^{g_\theta}(T,T)=\frac{m}{2}g_\theta(T,T),\qquad 
Ric^{g_\theta}(T,X^*)=0\qquad \mbox{and}\\[3mm]
Ric^{h}(X,Y)=Ric^{g_\theta}(X^*,Y^*)+\frac{1}{2}g_\theta(X^*,Y^*)\end{array}\] 
for all $X,Y\in TN$.

Using the above result for the Ricci tensor of $g_\theta$ 
with respect to the Webster-Ricci curvature, we obtain 
\[Ric^{h}(X,Y)=iRic^W(X^*,JY^*)\]
for all $X,Y\in 
TN$.
Basically, this result says that the Ricci-Webster curvature
of a (TSPH)-structure is the Ricci curvature of the base space of the natural 
submersion.

%%%%%%%%%%%%%%%%%%%%%%%%%%%%%%%%
\section{Description and construction of pseudo-Hermitian Einstein 
spaces}%%%
%%%%%%%%%%%%%%%%%%%%%%%%%%%%%%%%
\label{ab6}

We explain here an explicit constructions of pseudo-Hermitian Einstein 
spaces with arbitrary Webster scalar curvature. 
We also show that locally this construction principle generates
all pseudo-Hermitian Einstein structures. So we gain 
a locally complete description. The ideas in this section suggest 
that the construction can be extended to conformal K{\"a}hler 
geometry, in general. We aim to discuss this approach somewhere else.  

Let $(M,H,J,\theta)$ be a pseudo-Hermitian Einstein space with
arbitrary signature $(p,q)$, i.e. it holds
\[Ric^W=-i\frac{\ scal^W}{m}\cdot d\theta\qquad\mbox{and}\qquad
Tor^W(T,X)=0\]
for all $X\in\Gamma(H)$. Moreover, we assume for simplicity that $\theta$
generates globally a smooth Riemannian submersion
\[\pi:(M,g_\theta)\to (N,h)\ .\] 
With the relation for the Ricci tensors
from the end of the last section we obtain
\[\pi^*Ric^{h}=\frac{\ scal^W}{m} 
d\theta(\cdot,J\cdot)=\frac{\ scal^W}{m}\pi^*h\ .\]
This shows that the base space of the natural submersion to the
(TSPH)-structure $\theta$ is a K{\"a}hler-Einstein space of scalar 
curvature \[scal^{h}=2\cdot scal^W\ .\]
We conclude  that a pseudo-Hermitian Einstein space $(M,H,J,\theta)$
of dimension $n=2m+1$ 
determines uniquely (at least locally) a K{\"a}hler-Einstein 
manifold
$(N,h,J)$ of dimension $2m$ and signature $(2p,2q)$. 

We want to show now that there is a construction which assigns
to any K{\"a}hler-Einstein metric (with signature $(2p,2q)$) a
uniquely determined pseudo-Hermitian structure (which is then Einstein).
The construction itself is only unique up to gauge transformations.
However, it is easy to check that the resulting pseudo-Hermitian structures
to any gauge
are isomorphic.
To start with, 
let $(N^{2m},h,J)$ be a K{\"a}hler-Einstein space of dimension $2m$ with 
$scal^h>0$ and let 
$P(N)$ be the $\U(n)$-reduction of
the orthonormal frame bundle to $(N,h)$. Then it is
\[\mathcal{S}_{ac}(N):=P(N)\times_{det}S^1\]
the principal $S^1$-fibre bundle over $N$, which is associated to the
anti-canonical complex line bundle $\mathcal{O}(-1)$ of the K{\"a}hler 
manifold
$(N,h,J)$. The Levi-Civita connection to $h$ induces 
a connection form $\rho_{ac}$ on the anti-canonical $S^1$-bundle 
$\mathcal{S}_{ac}(N)$ with values in $i\RR$. For its curvature we have
\[\Omega^{\rho_{ac}}(\pi^{-1}_{\mathcal{S}_{ac}(N)*}X,
\pi^{-1}_{\mathcal{S}_{ac}(N)*}Y)
=iRic^h(X,JY),\qquad 
X,Y\in TN\ .\] 
At first, we see from this formula that the horizontal spaces of
$(\mathcal{S}_{ac}(N),\rho_{ac})$ generate a contact distribution $H$ of 
codimension
$1$ in $T\mathcal{S}_{ac}(N)$ and the 
horizontal lift of the complex structure
$J$ to $H$ produces a non-degenerate CR-structure $(H,J)$ on 
$\mathcal{S}_{ac}(N)$. This CR-structure is integrable as can be seen from 
the relation
\[\Omega^{\rho_{ac}}(X^*,JY^*)
+\Omega^{\rho_{ac}}(JX^*,Y^*)=0\]
for all $X,Y\in TN$ and the fact that the Nijenhuis tensor $N(X^*,Y^*)$ is 
the 
horizontal lift of 
\[J([JX,Y]+[X,JY])-[JX,JY]+[X,Y]=0\ .\]
Secondly, 
we see that
\[\theta:=i\frac{2m}{\ scal^h}\rho_{ac}\]
is a pseudo-Hermitian structure on $M:=\mathcal{S}_{ac}(N)$
furnished with the CR-structure $(H,J)$. 
The Reeb vector field $T$ on the 
pseudo-Hermitian space  
\[(\mathcal{S}_{ac}(N),H,J,\theta)\]
is vertical along the fibres (in fact, it is a fundamental vector field
generated by the right action) and by 
construction of $(H,J)$ transversally symmetric.  
Since $d\theta=\pi^*_{\mathcal{S}_{ac}(N)}h(J\cdot,\cdot)$ on $H$,  the 
base 
space of the 
corresponding
submersion is again the K{\"a}hler-Einstein space $(N,h,J)$ that we 
started
with. For that reason, we know that the Webster-Ricci curvature to 
$\theta$
must be given by
\[iRic^W(X^*,JY^*)=Ric^{h}(X,Y),\qquad X,Y\in H\ 
.\]
Since  $h$ is Einstein, we can conclude that the pseudo-Hermitian 
space
\[(\mathcal{S}_{ac}(N),H,J,\theta)\] 
is Einstein as well with Webster-Ricci curvature
\[Ric^W=-i\frac{\ scal^h}{2m}\cdot d\theta\ .\]

As mentioned before, for the inverse construction on the K{\"a}hler-Einstein space 
$(N,h,J)$, the choice of $\theta=i\frac{2m}{\ scal^h}\rho_{ac}$ as 
pseudo-Hermitian $1$-form is not unique.
One might replace $\theta$ by $\hat{\theta}:=\theta+df$ for some 
smooth function
$f$ on $\mathcal{S}_{ac}(N)$ with $df\neq -\theta$. 
The latter condition ensures that $\hat{\theta}$ is `transversal`,
which makes it possible to lift the complex structure to the kernel of
$\hat{\theta}$. 
We obtain again
a pseudo-Hermitian Einstein structure to $\mathcal{S}_{ac}(N)$ with
induced CR-structure.
It is straightforward to see that there is a diffeomorphism (gauge 
transformation)
on $\mathcal{S}_{ac}(N)$, which transforms $\theta+df$ into $\theta$, i.e. 
there is an isomorphism of pseudo-Hermitian structures. Since (locally) 
$\theta+df$ is 
the most general choice of a 'transversal' $1$-form whose exterior 
differential is the lift of $h(J\cdot,\cdot)$ on $N$, 
we
know that our gauged construction exhausts locally all pseudo-Hermitian 
Einstein
structures with positive Webster scalar curvature.

For the case of negative Webster scalar curvature note that if
$(M,H,J,\theta)$ with signature $(p,q)$ has positive Webster scalar 
curvature $scal^W>0$
then $(M,H,J,-\theta)$ has negative Webster scalar curvature 
$-scal^W$ and the base 
space of the natural submersion is $(N,-h,J)$, which is
K{\"a}hler-Einstein with reversed complex signature $(q,p)$. 
So if $(N,h,J)$ has $scal^h<0$ then $(N,-h,J)$ has $scal^{-h}>0$
and $\theta=i\frac{2m}{scal^{-h}}\rho_{ac}$ has positive Webster scalar 
curvature.
Hence, the pseudo-Hermitian form $i\frac{2m}{scal^h}\rho_{ac}$ has 
negative 
Webster scalar 
curvature. 
We conclude that
any pseudo-Hermitian Einstein structure  with $scal^W\neq 0$ can be 
realised 
locally
on $(\mathcal{S}_{ac}(N),H,J)$ over a K{\"a}hler-Einstein space with
$scal^h\neq 0$ by 
$\theta=i\frac{2m}{\ 
scal^h}\rho_{ac}$.  

As we have seen above a Webster-Ricci flat pseudo-Hermitian space 
$(M,H,J,\theta)$ gives rise to a Ricci-flat K{\"a}hler space. Again 
we aim to find an inverse construction. So let $(N,h,J)$ be a 
Ricci-flat K{\"a}hler space furnished with a $1$-form $\gamma$ 
such that $d\gamma=h(\cdot,J\cdot)$, i.e. $\omega:=d\gamma$ is
the K{\"a}hler form. 
The  $S^1$-principal fibre bundle $\mathcal{S}_{ac}(N)$ has a Levi-Civita 
connection form $\rho_{ac}$ with values in $i\RR$ which is flat, i.e. 
$d\rho_{ac}=0$. We set
\[\theta:=i\rho_{ac}-\pi^*\gamma\]
on $\mathcal{S}_{ac}(N)$. Obviously, it holds
\[d\theta=
-\pi^*\omega\ ,\] i.e. $\theta$ is a contact form
on $\mathcal{S}_{ac}(N)$ and the distribution $H$ in 
$T\mathcal{S}_{ac}(N)$,
which is given by $\theta|_H\equiv 0$ is contact as well. 
By definition, the distribution $H$ is transversal to the vertical 
direction of the 
fibre. For that reason we can lift $J$ to $H$. Again, the CR-structure
$(H,J)$ on $M:=\mathcal{S}_{ac}(N)$ is integrable.
Moreover, 
$\theta$ is a pseudo-Hermitian structure on $(\mathcal{S}_{ac}(N),H,J)$.
As the construction is done, it 
is clear that locally around every point of $(\mathcal{S}_{ac}(N),g_\theta)$
the base of the natural Riemannian submersion is a subset of the 
Ricci-flat space $(N,h,J)$. We conclude that
\[iRic^W(X,Y)=Ric^{h}(X,Y)=0\]
for all $X,Y\in H$, i.e. the pseudo-Hermitian space
\[(\mathcal{S}_{ac}(N),H,J,\theta)\] over a 
Ricci-flat 
K{\"a}hler space $(N,h,J)$ with K{\"a}hler form $d\gamma$, 
where $\theta=i\rho_{ac}-\pi^*\gamma$, is Webster-Ricci flat. 

In the Webster-Ricci flat construction, the pseudo-Hermitian form
$\theta$ can be replaced by $\hat{\theta}=i\rho_{ac}-\pi^*\gamma+df$, 
where 
$f$ is some smooth function on $\mathcal{S}_{ac}(N)$ with $df\neq 
-i\rho_{ac}$. 
This is the most
general 'transversal' $1$-form on $\mathcal{S}_{ac}(N)$ with 
$d\hat{\theta}=-\pi^*\omega$.  
However,
again one can see that  $\theta$ and 
$\hat{\theta}=\theta+df$ 
are gauge equivalent on $\mathcal{S}_{ac}(N)$, i.e. they are isomorphic
as pseudo-Hermitian structures. We conclude that with our construction of a particular gauge
we found (locally) the most general 
form of a Webster-Ricci flat pseudo-Hermitian  space. 
We summarise our results.

\begin{THEO} \label{TH1} Let $(N,h,J)$ be a K{\"a}hler-Einstein space of 
dimension $2m$ and signature $(2p,2q)$ with scalar curvature $scal^h$. 
\begin{enumerate} 
\item If $scal^h\neq 0$ then
the anti-canonical $S^1$-principal bundle 
\[\mathcal{S}_{ac}(N)=P(N)\times_{det}S^1\] with canonically
induced CR-structure $(H,J)$ and connection $1$-form
\[\theta:=i\frac{2m}{\ scal^h}\rho_{ac}\ ,\] where $\rho_{ac}$ is the 
Levi-Civita 
connection to $h$, 
is a pseudo-Hermitian Einstein space with 
$scal^W=\frac{1}{2}scal^h\neq 0$.
\item
If $scal^h=0$ and the K{\"a}hler form is
$\omega=d\gamma$ for some $1$-form $\gamma$ on $N$ 
then $(\mathcal{S}_{ac}(N),H,J)$ with 
pseudo-Hermitian structure $\theta=i\rho_{ac}-\pi^*\gamma$ 
is Webster-Ricci flat. 
\end{enumerate}
Locally, any pseudo-Hermitian Einstein space 
$(M,H,J,\theta)$
is isomorphic to one of these two models depending
on the Webster scalar curvature $scal^W$.
\end{THEO} 

We remark here that for the case $scal^h\neq 0$      
we could have chosen the gauge 
$\theta=i\rho_{ac}+\pi^*\eta-\pi^*\gamma$, where $d\gamma$ is the
K{\"a}hler form and $d\eta$ the Ricci form. This would enable
us to treat the two cases of Theorem \ref{TH1} as one case.
However, for the following discussion of the 
corresponding Fefferman spaces we find the chosen 
gauge of Theorem \ref{TH1} more convenient.

%%%%%%%%%%%%%%%%%%%%%%%%%%%%%%%%
\section{The Fefferman metric to a pseudo-Hermitian structure}%%%
%%%%%%%%%%%%%%%%%%%%%%%%%%%%%%%%
\label{ab7}
We briefly explain here the construction of the Fefferman space 
which belongs to any pseudo-Hermitian space (cf. 
\cite{Fef76}, \cite{Spa85}, \cite{Lee86}, \cite{Bau99}).

Let $(M,H,J,\theta)$ be a pseudo-Hermitian space of dimension 
$n=2m+1$ and signature $(p,q)$. We denote by $(F,\pi_M,M)$ the 
canonical
$S^1$-principal fibre bundle of the CR-manifold $(M,H,J)$. 
The Webster connection 
on $F$ is denoted by $A^W$. In general, it holds
\[\Omega^W=dA^W=-\pi_M^*Ric^W \ .\]
Furthermore, we define
\[A_\theta:=A^W-\frac{i}{2(m+1)}scal^W\pi^*_M\theta\ .\]
The Fefferman metric on $F$ belonging to $\theta$ on $(M,H,J)$ is defined 
as
\[ f_\theta:= \pi_M^*L_\theta-i\frac{4}{m+2}\pi_M^*\theta\circ A_\theta\ 
.\]
The signature of this metric is $(2p+1,2q+1)$. The $1$-forms $\pi_M^*\theta$ 
and $A_\theta$ are both lightlike
with respect to the Fefferman metric $f_\theta$.
In particular, if 
$(M,H,J,\theta)$ is strictly pseudoconvex, the signature of 
$(F,f_\theta)$ is Lorentzian.

Let $P$ denote the fundamental vector field in vertical direction on $F$ 
generated by
the element $\frac{m+2}{2}i\in i\RR$, i.e. $A_\theta(P)=\frac{m+2}{2}i$. 
Moreover, in this section we denote by
$X^*$ the horizontal lift with respect to 
$A_\theta$ of a vector field $X$ in $\Gamma(H)$ on $M$. With $T^*$ we 
denote 
the horizontal lift of the Reeb vector field $T$. With our definitions 
it is
\[f_\theta(P,T^*)=1\ .\]
%The vector field $\frac{1}{\sqrt{2}}(N-T^*)$ has length $-1$ with respect 
%to $f_\theta$ and the vector field $\frac{1}{\sqrt{2}}(N+T^*)$
%has length $1$. The vector field $N$ is a lightlike Killing vector field
%on  $(F,f_\theta)$. 
We have the following formulas for commutators
and covariant derivatives with respect to $f_\theta$ (cf. \cite{Bau99}):
\[\begin{array}{l}
[X^*,P]=[T^*,P]=0\ ,\\[2mm]
Vert_\pi [X^*,Y^*]=i\frac{2}{m+2}\Omega^{A_\theta}(X^*,Y^*)\cdot P\ 
,\\[2mm]
Horiz_\pi [X^*,Y^*]=[X,Y]^*\ ,\\[2mm]   
[T^*,X^*]=[T,X]^*+i\frac{2}{m+2}\Omega^{A_\theta}(T^*,X^*)\cdot P\ 
,\\[2mm]
[X^*,Y^*]=\pi_H[X,Y]^*-d\theta(X,Y)\cdot 
T^*+i\frac{2}{m+2}\Omega^{A_\theta}(X^*,Y^*)\cdot P\ .
\end{array}\]
Furthermore,
\[\begin{array}{l}
f_\theta(\nabla^{f_\theta}_{X^*}Y^*,Z^*)=L_\theta(\nabla^W_XY,Z)\ ,\\[2mm]
f_\theta(\nabla^{f_\theta}_PY^*,Z^*)=\frac{1}{2}d\theta(Y,Z)\ ,\\[2mm]
f_\theta(\nabla^{f_\theta}_{T^*}Y^*,Z^*)=\frac{1}{2}\left(
L_\theta([T,Y],Z)-L_\theta([T,Z],Y)-i\frac{2}{m+2}
\Omega^{A_\theta}(Y^*,Z^*)\right)\ ,\\[2mm] 
f_\theta(\nabla^{f_\theta}_{X^*}Y^*,P)=-\frac{1}{2}d\theta(X,Y)\ ,\\[2mm]
f_\theta(\nabla^{f_\theta}_{X^*}Y^*,T^*)=\frac{1}{2}\left(
L_\theta([T,X],Y)+L_\theta([T,Y],X)+i\frac{2}{m+2}
\Omega^{A_\theta}(X^*,Y^*)\right)\ ,\\[2mm]
f_\theta(\nabla^{f_\theta}_{T^*}T^*,Z^*)=-i\frac{2}{m+2}\Omega^{A_\theta}
(T^*,Z^*)\ 
,\\[2mm]
f_\theta(\nabla^{f_\theta}P,T^*)=f_\theta(\nabla^{f_\theta}T^*,T^*)
=f_\theta(\nabla^{f_\theta} P,P)=0\ ,\\[2mm]
f_\theta(\nabla^{f_\theta}_PP,Z^*)=f_\theta(\nabla^{f_\theta}_PT^*,Z^*)
=f_\theta(\nabla^{f_\theta}_{T^*}P,Z^*)=0\end{array}\]
for all $X,Y,Z\in \Gamma(H)$.

%%%%%%%%%%%%%%%%%%%%%%%%%%%%%%%%%%%%%%%%%%%%%%%%%%%%%%%
\section{Einstein-Fefferman spaces}
\label{ab8}
%%%%%%%%%%%%%%%%%%%%%%%%%%%%%%%%%%%%%%%%%%%%%%%%%%%%%%%%%

We discuss here the Fefferman metric of a pseudo-Hermitian
Einstein space. In particular, we give an explicit local 
construction for an Einstein metric in the conformal class of any 
Fefferman metric coming from a pseudo-Hermitian Einstein space.

Let $(M,H,J,\theta)$ be a pseudo-Hermitian Einstein
space. We know already that every 
such pseudo-Hermitian Einstein space is constructed at least locally from 
a K{\"a}hler-Einstein space. We assume here for simplicity that
$M=\mathcal{S}_{ac}(N)$ is the total space of the anti-canonical 
$S^1$-principal fibre bundle over a K{\"a}hler-Einstein 
space $(N,h,J)$
furnished with the naturally induced CR-structure $(H,J)$ and 
pseudo-Hermitian form $\theta$ as described in 
Theorem \ref{TH1}:
\[\pi_N^{ac}: (\mathcal{S}_{ac}(N),H,J,\theta)\ \to\ (N,h,J)\ 
.\]
Moreover, we denote by 
\[\mathcal{S}_{c}(N):=P(N)\times_{det^{-1}}S^1\]
the canonical $S^1$-principal fibre bundle over $(N,h,J)$ which
is furnished with the Levi-Civita connection denoted by $\rho_c$.
Now let $F$ be the total space of the canonical $S^1$-fibre bundle
over the CR-manifold $(\mathcal{S}_{ac}(N),H,J)$. We denote by
$\pi$ the projection of $F$ to $N$:
\[\pi:F\to N\ .\]
Obviously, 
the lift of 
$\mathcal{S}_{c}(N)$
along the anti-canonical projection $\pi_N^{ac}$ is isomorphic to $F$. 
This 
shows that we can 
understand the 
total space
$F$ as a torus bundle over $(N,h,J)$: 
\[F=P(N)\times_{(det,det^{-1})}S^1\times S^1\ .\]
On $F$ we already introduced the $1$-forms $\pi^*_{ac}\theta$
and $A^W$ resp. $A_\theta$. In our situation here, where the Fefferman 
construction is based on a K{\"a}hler-Einstein space, we can 
express the two latter $1$-forms on $F$ by:
\[A^W=\pi_{c}^*\rho_c\qquad\mbox{and}\qquad A_\theta=\pi_c^*\rho_c-
\frac{\ 
i\cdot scal^W}{2(m+1)}\pi_{ac}^*\theta
%=\pi_c^*\rho_c+\frac{m}{2(m+1)}\pi_{ac}^*\rho_{ac}
\ ,\]
where $\pi_c:F\to \mathcal{S}_c(N)$ and 
$\pi_{ac}:F\to \mathcal{S}_{ac}(N)$ are the natural projections.
  
As it was defined in the previous section, the
Fefferman metric $f_\theta$ to the pseudo-Hermitian Einstein space
$(\mathcal{S}_{ac}(N),H,J,\theta)$ lives on the total space $F$.
We can express now the Fefferman metric on $F$ in the pseudo-Hermitian 
Einstein
case by
\[f_\theta=\pi^*h-i\frac{4}{m+2}
\pi^*_{ac}\theta\circ
\left(\pi_c^*\rho_c-i\frac{\ scal^W}{2(m+1)}
\pi_{ac}^*\theta\right)\ .\] 
%Thereby, we switched the notation for the Fefferman metric from  
%$f_\theta$ to $f_h$ to 
%indicate that this metric is basically derived solely from the underlying  
%K{\"a}hler-Einstein space $(N,h,J)$.
Notice that the metric $f_\theta$ is uniquely derived from $(N,h,J)$ if we 
assume
$\theta$ to be given in the gauge as described in Theorem \ref{TH1}.
For simplicity we will in the following omit
the subscripts of the projections: $\pi=\pi_c$ and $\pi=\pi_{ac}$. 
It will be clear from the context which projection is meant. 
\begin{DEF} Let $(N,h,J)$ be a K{\"a}hler-Einstein 
space and let 
\[F=P(N)\times_{(det,det^{-1})}S^1\times 
S^1\] be the `canonical-anti-canonical' torus bundle over $N$.
Then we denote by  $f_h:=f_\theta$ (where $\theta$ is the gauge given as 
in 
Theorem \ref{TH1} depending on $scal^h$)  
the Fefferman metric on $F$ 
which belongs to the pseudo-Hermitian Einstein space 
$(\mathcal{S}_{ac}(N),H,J,\theta)$. We call $f_h$ the Fefferman metric
of the K{\"a}hler-Einstein  space $(N,h,J)$.
\end{DEF}

In general, it is 
\[d(\pi^*\rho_c)+d(\pi^*\rho_{ac})=\pi^*Ric^W-\pi^*Ric^W=0\ ,\] 
i.e. the $1$-form $\pi^*\rho_c+\pi^*\rho_{ac}$ is closed on $F$. 
In fact, we will see below that this $1$-form is parallel in the 
Einstein case.
If $scal^h=0$ we calculate the Fefferman metric of $N$ from the
above expression to
\[f_h=
\pi^*h-i\frac{4}{m+2}(i\pi^*\rho_{ac}-\pi^*\gamma)\circ\pi^*\rho_c\
,\]
where $d\gamma=\omega$ is the K{\"a}hler form. In case that $scal^h\neq 0$
an orthogonal $1$-form to $\pi^*\rho_c+\pi^*\rho_{ac}$ is given by 
$\pi^*\rho_c-\frac{1}{m+1}\pi^*\rho_{ac}$ and we can write  
the Fefferman metric as
\[f_h=\pi^*h+\frac{4m(m+1)}{(m+2)^2\cdot scal^h}
\cdot\left(\ (\pi^*\rho_c+\pi^*\rho_{ac})^2
-(\pi^*\rho_c-\frac{1}{m+1}\pi^*\rho_{ac})^2\ 
\right)\ .\]

We want to  calculate the Ricci tensor of $f_h$. 
As before, let $P$ be 
the vertical vector field on $F$ along the 'Fefferman' fibering with 
\[A^W(P)=\pi^*\rho_c(P)=\frac{m+2}{2}i\qquad\mbox{and}\qquad
\pi^*\rho_{ac}(P)=0\] 
and let $T^*$ be the vertical vector field along the anti-canonical
fibering with
\[\pi^*\theta(T^*)=1\qquad\mbox{and}\qquad 
A_\theta(T^*)=0\ .\]
Furthermore, let 
\[(e_i)_{i=1,\ldots,2m}\]
denote a local orthonormal basis on $(N,h,J)$ and let $e_i^*$ 
be the horiziontal lifts of $e_i$ to $F$ with respect to $\theta$ and 
then $A_\theta$, 
i.e. 
\[ \pi_*(e_i^*)=e_i\qquad\mbox{and}\qquad 
\pi^*\theta(e_i^*)=\pi^*\rho_c(e_i^*)=0\ .\] 
Then we 
have 
\[ [T^*,e_i^*]=[P,e_i^*]=[P,T^*]=0 \qquad \mbox{for\ all}\ \ i\in 
1,\ldots,2m\]
on $F$. 
We will work in the
following always with a local
basis on $F$
of the form \[(e_i^*,T^*,P)\ .\]

Now, since $\theta$ is Einstein, we observe that
\[dA_\theta=\Omega^{A_\theta}=-\pi_{ac}^*Ric^W
-i\frac{ 
scal^W}{2(m+1)}\pi_{ac}^*d\theta=i\frac{(m+2)\cdot 
scal^h\ }{4m(m+1)}\cdot
\pi_{ac}^*d\theta\ 
.\]
We use this and the formulas from the last section (cf. \cite{Bau99}) to 
obtain the
covariant derivatives for a local basis $(e_i^*,T^*,P)$ on $(F,f_h)$.
It is
\[\begin{array}{l}
\nabla^{f_h}_{e_i^*}e_j^*=(\nabla^W_{e_i}e_j)^*-\frac{1}{2}
\pi^*d\theta(e_i^*,e_j^*)T^*-\frac{1}{2}S^W\pi_{ac}^*d\theta(e_i^*,e_j^*)
P\\[3mm]
\nabla^{f_h}_{T^*}e_i^*=\nabla^{f_h}_{e_i^*}T^*
=\frac{1}{2}S^W(Je_i)^* 
\\[3mm]
\nabla^{f_h}_Pe_i^*=\nabla^{f_h}_{e_i^*}P=\frac{1}{2}(Je_i)^*\\[3mm]
\nabla^{f_h}_PT^*=\nabla^{f_h}_{T^*}P=\nabla^{f_h}_{T^*}T^*=
\nabla^{f_h}_PP=0\ ,
\end{array}\]
whereby we set \[S^W:=\frac{scal^h}{2m(m+1)}\ .\]
It follows immediately that 
\[f_h(\nabla^{f_h}_{A}T^*,B)=
-f_h(\nabla^{f_h}_{B}T^*,A)\]
for all $A,B\in \Gamma(TF)$, i.e. $T^*$ is a Killing vector field on 
$(F,f_h)$. (In general, for any pseudo-Hermitian space the 
horizontal lift of the Reeb vector field $T$ is a Killing vector
on the Fefferman space if and only if $T$ is a transversal symmetry
and $\Omega^{A_\theta}(T^*,\cdot)=0$.)

From the formulas for the covariant derivative we see that
the vertical vector field
\[T^*-S^WP\]
is parallel. The dual of this vector field with respect to 
the Fefferman metric $f_h$ 
is equal to $-i\frac{2}{m+2}\cdot(\pi^*\rho_c+\pi^*\rho_ac)$,
which is a parallel $1$-form.
For the Riemannian curvature tensor of $f_h$ we find
\[\begin{array}{l}
R^{f_h}(e_i^*,e_j^*)e_j^*=(R^{\nabla^W}(e_i,e_j)e_j)^*+
\frac{3}{2}S^Wd\theta(e_i,e_j) (Je_j)^*\ ,\\[3mm]
R^{f_h}(e_i^*,P)T^*=\frac{1}{4}S^W\cdot e_i^*\ ,\\[3mm]
R^{f_h}(T^*,e_j^*)e_j^*=\frac{1}{4}S^W\cdot (T^*+S^W\cdot P)\ 
,\\[3mm]
R^{f_h}(P,e_j^*)e_j^*=\frac{1}{4}(T^*+S^W\cdot P)\ ,\\[3mm]
R^{f_h}(P,T^*)=0\ .\end{array}\]
Then we obtain for the Ricci tensor 
\[\begin{array}{l} 
Ric^{f_h}(e_i^*,e_j^*)=iRic^W(e_i,Je_j)-S^W
f_h(e_i^*,e_j^*)\ ,\\[3mm]
Ric^{f_h}(T,e_i^*)=Ric^{f_h}(P,e_i^*)=0\ ,\\[3mm]
Ric^{f_h}(T^*,T^*)=\frac{m}{2}(S^W)^2\ ,\\[3mm]
Ric^{f_h}(T^*,P)=\frac{m}{2}S^W\ ,\\[3mm]
Ric^{f_h}(P,P)=\frac{m}{2}\ ,\end{array}\]
i.e. the Ricci tensor of $f_h$ takes the form
\begin{eqnarray*}Ric^{f_h}&=\quad&i\pi^*Ric^W(\cdot,J\cdot)
-S^W\pi^*f_h\\[2mm]&\quad +&\frac{m}{2}\left(
(S^W)^2
\pi^*\theta\circ
\pi^*\theta-\frac{4}{(m+2)^2}A_\theta\circ
A_\theta-i\frac{4}{m+2}S^WA_\theta\circ\pi^*\theta\right)\\[3mm]
&=\quad&\frac{scal^h}{2(m+1)}f_h-\frac{2m}{(m+2)^2}\left(A_\theta
-\frac{\ i(m+2)\cdot scal^h\ }{4m(m+1)}\pi^*\theta\right)^2\ . 
\end{eqnarray*}
In particular, if $scal^h=0$ on $N$ then
\[Ric^{f_h}=-\frac{2m}{(m+2)^2}(\pi_c^*\rho_c)^2\]
and if $scal^h\neq 0$ then
\[Ric^{f_h}=
\frac{scal^h}{2(m+1)}f_h-\frac{2m}{(m+2)^2}(\pi^*\rho_c+\pi^*\rho_{ac})^2\ 
.\]

The calculations show that the Fefferman metric $f_h$ to a 
K{\"a}hler-Einstein space
$(N,h,J)$  
is never Einstein (cf. \cite{Lee86}). For example, a Webster-Ricci flat 
pseudo-Hermitian 
space gives rise to a Fefferman metric $f_h$ with totally isotropic 
Ricci tensor. However, the 
Einstein condition should not be expected for the Fefferman metric.
Instead, we will show now that the Fefferman metric to any 
pseudo-Hermitian 
Einstein space (resp. K{\"a}hler-Einstein space) is (locally) conformally 
Einstein,
i.e. there is locally a conformally rescaled metric $\tilde{f}_h$
to $f_h$ 
which is Einstein. For the calculation of the  conformal Einstein scale
we introduce the  coordinate function $t$ on the torus fibre bundle $F$
by  
\[\begin{array}{ll}dt=i\pi^*\rho_c\qquad\qquad& \mbox{when}\quad 
scal^h=0\quad \mbox{and}\\[2mm]
dt=i\pi^*\rho_c+i\pi^*\rho_{ac} \qquad\qquad& \mbox{when}\quad scal^h\neq 
0\ .
\end{array}\]

First, we consider the Webster-Ricci flat case. For simplicity we 
assume  a 
submersion 
\[\pi:(\mathcal{S}_{ac}(N),H,J,\theta)\to (N,h,J) \ ,\]
where $\theta=i\pi^*\rho-\pi^*\gamma$ and $\omega=d\gamma$
is the K{\"a}hler form.
Let $f_h$ be the Fefferman metric on $F$ over $N$ and let 
$\tilde{f}_h=e^{2\phi}f_h$ be a conformally rescaled metric
with real function $\phi$ on $F$. For the Ricci
tensor of $\tilde{f}_h$ we find by using standard formulas 
\[Ric^{\tilde{f}_h}-Ric^{f_h}=-2m(Hess(\phi)-d\phi\circ d\phi)+
(-\Delta \phi-2m\|d\phi\|^2)f_h\ ,\]
where $\Delta$ denotes the Laplacian with respect to $f_h$. We denote
the correction term on the right hand side of this formula by $C_\phi$.
This 
is a 
symmetric 
$2$-tensor.
Let $\phi(t)$ be a function on $F$ which depends only
on the coordinate $t$ in direction of  the canonical $S^1$-fibering.
Then the only non-trivial component of $C_\phi$ is
\[C_\phi(P,P)=-2m\left(\ PP(\phi)-P(\phi)^2\ \right)\ . \] 
This shows that for any function $\phi(t)$ on 
$F$, which satisfies
the ODE
\[\partial_{t}\partial_{t}\phi-(\partial_{t}\phi)^2=\frac{1}{(m+2)^2}\ 
,\]
the Ricci tensor $Ric^{\tilde{f}_h}$ of $\tilde{f}_h=e^{2\phi}f_h$
vanishes. The most general solution of this ODE is
\[\phi=c_1-ln\left(cos\left(\frac{t}{m+2}+c_2\right)\right) ,\]
where $c_1, c_2$ are constants. We choose here
$\phi=-ln(cos(\frac{t}{m+2}))$, which then gives as conformal rescaling
factor
\[e^{2\phi}= cos^{-2}\left(\frac{t}{m+2}\right)\ .\]
Then the conformally changed Fefferman metric (in short: conformally 
Fefferman 
metric) 
\[\tilde{f}_h=cos^{-2}(t/(m+2))\cdot    
\left(\pi^*h-i\frac{4}{m+2}(i\pi^*\rho_{ac}
-\pi^*\gamma)\circ\pi^*\rho_c\right) 
\]
is Ricci-flat on an open subset in $F$ around the hypersurface 
given by $\{t=0\}$. 
Obviously, a global conformal Einstein scale for $f_h$ on $F$ does not 
exist. Everywhere locally it does exists.

We assume now that  $(N,h,J)$ is K{\"a}hler-Einstein 
with $scal^h\neq 0$. Then we find with respect to a conformal 
scaling function $\phi(t)$, which depends only on the coordinate $t$
with $dt=i(\pi^*\rho_c+\pi^*\rho_{ac})$,

\begin{eqnarray*}
Ric^{\tilde{f}_h}-Ric^{f_h}=\ C_\phi 
&=&\quad
2m(\partial_t\partial_t\phi-(\partial_t\phi)^2)(d\rho_c+d\rho_{ac})^2\\[2mm]
&&+\ \frac{(m+2)^2\cdot 
scal^h}{4m(m+1)}(\partial_t\partial_t\phi-(\partial_t\phi)^2)f_h\ .
\end{eqnarray*}
Again, if we choose $\phi=-ln(cos(\frac{t}{m+2}))$, 
the metric $\tilde{f}_h=e^{2\phi}f_h$ is Einstein. In fact, we obtain
with this function $\phi$ for the Ricci tensor of the rescaled metric
\[Ric^{\tilde{f}_h}=\frac{(2m+1)\cdot scal^h}{4m(m+1)}f_h\]
and the scalar curvature is 
$scal^{\tilde{f}_h}=\frac{2m+1}{2m}\cdot scal^h$.
\begin{THEO} \label{TH2} Let $(N,h,J)$ be a K{\"a}hler-Einstein space of
dimension $2m$ and signature $(2p,2q)$ with scalar curvature $scal^h$.
\begin{enumerate}
\item If $scal^h=0$ and $\omega=d\gamma$ for some $1$-form $\gamma$ on $N$ 
then the metric
\[\tilde{f}_h=cos^{-2}(t)\cdot\left(\ \pi^*h+4dt\circ(\pi^*\gamma+ds)\ 
\right)\\[5mm] \]
on $N\times\{\ (s,t)\; :\ -\frac{\pi}{2}< t< \frac{\pi}{2}\ \}\subset
N\times \RR^2$ (with natural projection $\pi$ onto $N$) 
is conformally Fefferman and Ricci-flat with signature $(2p+1,2q+1)$. 
\item
If $scal^h\neq 0$ then the metric
\[\tilde{f}_h=cos^{-2}(t)\cdot \left(\ \pi^*h-\frac{4m(m+1)}{\ scal^h\ }
\cdot(\ dt^2+\frac{\rho_{ac}^2}{(m+1)^2})\right)\\[3mm]\]
on $\mathcal{S}_{ac}(N)\times]-\frac{\pi}{2},\frac{\pi}{2}[$,
where $(\mathcal{S}_{ac}(N),\pi,N)$ is the anti-canonical $S^1$-bundle 
over
$N$ with Levi-Civita connection $\rho_{ac}:T\mathcal{S}_{ac}(N)\to i\RR$,
is conformally Fefferman and Einstein with 
$scal^{\tilde{f}_h}=
\frac{2m+1}{2m}\cdot scal^h$ and signature $(2p+1,2q+1)$. 
\end{enumerate}
Every Fefferman metric, which is conformally Einstein, is locally 
conformally equivalent to a metric of the form $\tilde{f}_h$ as described
here in (1) resp. (2).
\end{THEO}

In Theorem \ref{TH2},
we simplified the expressions for the 
Fefferman metrics.
In the Ricci-flat case both the 
Levi-Civita connections $\rho_c$ and $\rho_{ac}$ are flat, i.e. the torus
bundle is globally a product and we parametrised the vertical directions 
by the coordinates $t,s$, where the coordinate $t$ is rescaled 
(compared with our notation before) by a factor $(m+1)$ .

For the case when $scal^h\neq 0$ we replaced the $1$-form 
$\pi^*\rho_c-\frac{1}{m+1}\pi^*\rho_{ac}$ on $F$ by 
$-\frac{m+2}{m+1}\cdot \rho_{ac}$ on $\mathcal{S}_{ac}(N)$.
This is possible,  since locally the canonical bundel $\mathcal{S}_c(N)$
and the anti-canonical bundle $\mathcal{S}_{ac}(N)$ can be identified
such that the Levi-Civita connection $\rho_c$ becomes $-\rho_{ac}$. It is 
useful to note here that the Fefferman metric $f_h=cos^2(t)\tilde{f}_h$ 
as presented in Theorem \ref{TH2} is the product of a 
real line with the metric
\[\pi^*h-\frac{4m(m+1)}{\ (m+1)^2\cdot scal^h\ }\cdot\rho_{ac}^2\ .\]
This is the well-known Einstein-Sasaki metric which is constructed from the 
K{\"a}hler-Einstein metric $h$.  

For the proof of Theorem \ref{TH2}, we remark that it does not follow yet
from our discussion that any conformally Einstein Fefferman metric
comes from a pseudo-Hermitian Einstein space. To see this point in the 
proof we can use the argument from tractor calculus, which says
that the parallel standard tractor on the Einstein-Fefferman space gives 
rise to a parallel standard tractor on the underlying CR-space. This 
implies that the CR-space admits a pseudo-Hermitian Einstein structure.

\section*{Acknowledgments}
I would like to thank A. Rod Gover for many helpful discussions 
and correspondence and bringing my attention to this topic. 
%and
%inviting him for a stay at the Univesity of Auckland. 
In fact, the idea 
for this work originated
during a visit to the University of Auckland in Auckland/New 
Zealand and is part of a broader project on CR- and conformal Einstein 
geometry. I am grateful to Rod for the invitation and I would 
like to thank for hospitality during my 
stay, which  was supported by ARGs Royal Society of New Zealand,
Marsden Grant no.$\backslash$ 02-UOA-108, and by the New Zealand 
Institute of 
Mathematics and its Applications.

%-----------------------------------------------------------------------

\end{sloppypar}

\begin{thebibliography}{1111111}
%-----------------------------------------------------------------------

 	
\bibitem[Bau99]{Bau99} H. Baum. {\it Lorentzian twistor spinors 
and CR-geometry}. Differential Geom. Appl. 11 (1999), no. 1, 69--96. 
 	
\bibitem[Fef76]{Fef76} C. Fefferman. {\it Monge-Ampere equations,
the Bergman kernel, and geometry of pseudoconvex domains}. Ann. Math. 103
(1976), 395-416.

\bibitem[Spa85]{Spa85} G.A.J. Sparling. {\it Twistor theory and the 
characterisation of Fefferman's conformal structures.} Preprint Univ.
Pittsburg, 1985.

\bibitem[Lee88]{Lee88} J. M. Lee. {\it Pseudo-Einstein structures on CR 
manifolds}.  Amer. J. Math.  110  (1988),  no. 1, 157--178.
 
\bibitem[Lee86]{Lee86} J. M. Lee. {\it The Fefferman metric and 
pseudo-Hermitian invariants}.  Trans. Amer. Math. Soc.  296  (1986),  no. 
1, 411--429. 
 	
									
\bibitem[Cap01]{Cap01} A. Cap. {\it Parabolic geometries, 
CR-tractors, and the Fefferman construction}. 8th International Conference 
on Differential Geometry and its Applications (Opava, 2001).  Differential 
Geom. Appl.  17  (2002),  no. 2-3, 123--138.
 
\bibitem[CG02]{CG02}
A. Cap, A. R. Gover. {\it Tractor calculi for 
parabolic geometries}.  Trans. Amer. Math. Soc.  354  (2002),  no. 4, 
1511--1548 (electronic). 


\bibitem[CS00]{CS00} A. Cap, H. Schichl. 
{\it Parabolic geometries and canonical Cartan connections}. 
Hokkaido Math. J. 
29 (2000), no. 3, 453--505. 

								 
\bibitem[ONe66]{ONe66} B. O'Neill. {\it The fundamental equations 
of a submersion}. Michigan Math. J. 13 (1966), 459--469. 

\bibitem[Gov04]{Gov04} A. R. Gover. {\it Almost conformally
Einstein manifolds and obstructions}. preprint 2004.

\bibitem[Lei05]{Lei05} F. Leitner. {\it Conformal Killing forms
with normalisation condition}. to appear in Rend.
Circ. Mat. Palermo (2) Suppl. No. 65 (2005).


\end{thebibliography}
\end{document}